\documentclass{amsart}
\newtheorem{thm}{Theorem}

\newcommand{\Rho}{\mathrm{P}}
\begin{document}
\title{Uniqueness of the stereographic embedding}
\author{Michael Eastwood}
\address{\hskip-\parindent
Mathematical Sciences Institute\\
Australian National University\newline
ACT 0200, Australia}
\email{meastwoo@member.ams.org}
\subjclass{Primary 53A30; Secondary 53C22.}
\begin{abstract}
The standard conformal compactification of Euclidean space is the round 
sphere. We use conformal geodesics to give an elementary proof that this is 
the only possible conformal compactification.
\end{abstract}
\renewcommand{\subjclassname}{\textup{2010} Mathematics Subject Classification}
\maketitle

\section{Introduction}\label{intro}
At the $34^{\mathrm{th}}$ Winter School on Geometry and Physics, held as always
in Srn\'{\i} in the Czech Republic, Charles Franc\`es gave a plenary series of
talks on `Conformal boundaries in pseudo-Riemannian geometry.' He discussed the
results contained in~\cite{F} and, in particular, the notions of {\em conformal
embedding\/}
$$(M,[g])\hookrightarrow^{\mathrm{open}}(N,[h])$$
and {\em conformal maximality\/} (where there is no non-trivial conformal 
embedding). The prototype conformal embedding is
$$({\mathbb{R}}^n,[\mbox{Euclidean metric}])\hookrightarrow
(S^n,[\mbox{round metric}])$$
where the image is the complement of a single point, which we shall write
as~$\infty$, and the mapping in the other direction
$S^n\setminus\{\infty\}\to{\mathbb{R}}^n$ is the standard stereographic
projection, which is well-known to be conformal. It is proved in \cite{F} that 
this is the only non-trivial conformal embedding of Euclidean space. In 
particular, it is the only possible {\em conformal compactification\/} 
of~${\mathbb{R}}^n$ as follows.
\begin{thm}\label{main_thm}
Suppose $n\geq 3$ and $(M,g)$ is a compact Riemannian $n$-manifold and
$\iota:{\mathbb{R}}^n\hookrightarrow M$ realises ${\mathbb{R}}^n$ as an open
subset of $M$ such that $\iota^*g$ is conformal to the standard Euclidean
metric on~${\mathbb{R}}^n$. Then $M$ must be the $n$-sphere~$S^n$, the metric
$g$ must be conformal to the standard round metric, and the embedding must be
inverse to the standard stereographic projection.
\end{thm}
Of course, this is the expected result and in \cite{F} two proofs are
presented. One proof \cite[Theorem~1.4]{F} is via a removable singularities
theorem based on the complement of ${\mathbb{R}}^n\hookrightarrow S^n$ having
small Hausdorff dimension. The other proof \cite[Proposition~2.3]{F} is based
on B.~Schmidt's notion of {\em conformal boundary\/} as further developed by
Franc\`es.

This article gives a proof of Theorem~\ref{main_thm} using {\em conformal
geodesics\/}. Apart from the elementary nature of this proof, it seems to be a
natural one and a good application of the theory conformal geodesics, which are
perhaps under-used in conformal differential geometry (unfortunately, as 
lamented in~\cite{T}, their behaviour is not so well understood in general).

Although Theorem~\ref{main_thm} also holds when $n=2$, the result is of a
different nature (global complex analysis rather than local to global
differential geometry). Also note that the corresponding theorem is false for
projective differential geometry: there are two distinct projective
compactifications of ${\mathbb{R}}^n$, namely
${\mathbb{R}}^n\hookrightarrow{\mathbb{RP}}_n$ as a standard affine
co\"ordinate patch and ${\mathbb{R}}^n\hookrightarrow S^n$ onto an open
hemisphere by inverse central projection.

This article starts with an exposition of conformal geodesics, which is more
than is strictly needed for the proof of Theorem~\ref{main_thm}. Hopefully,
this exposition will be useful for their further application in conformal
geometry.

I would like to thank Vladim\'{\i}r Sou\v{c}ek for his invitation to Charles
University and the Czech Winter School. I would also like to thank Charles
Franc\`es for an inspiring series of talks and for useful conversations.

\section{Conformal Geodesics}
Also known as {\em conformal circles\/}~\cite{BE}, these special curves may be
defined as follows. For any metric $g_{ab}$ in the conformal class, a curve
$\gamma\hookrightarrow M$ with parameterization $\tau:\gamma\to{\mathbb{R}}$ is
said to be a {\em conformal geodesic\/} if and only if
\begin{equation}\label{circles}
\partial A^a=3\frac{V.A}{V.V}A^a-\frac{3A.A}{2V.V}V^a+V.V\Rho_b{}^aV^b
-2\Rho_{bc}V^bV^cV^a,\end{equation}
where
\begin{itemize}
\item $V^a$ is tangent along $\gamma$ so that $V^a\nabla_a\tau=1$, the 
`velocity vector,'
\item $\partial=\partial/\partial\tau\equiv V^a\nabla_a$ acting on any tensor
field along $\gamma$,
\item $\nabla_a$ is the Levi-Civita connection for $g_{ab}$,
\item $A^a\equiv\partial V^a$, the `acceleration vector,'
\item $\Rho_{ab}$ is the conformal Schouten tensor 
$$\Rho_{ab}=\frac1{n-2}\left(R_{ab}-\frac{1}{2(n-1)}Rg_{ab}\right),$$
\item $R_{ab}$ is the Ricci tensor and $R$ the scalar curvature,
\item $X.Y$ is the inner product $X^aY_a$ of two fields $X^a$ and $Y^a$. 
\end{itemize}
The equation (\ref{circles}) is an ordinary differential equation for both the
curve $\gamma$ and its parameterization $\tau:\gamma\to{\mathbb{R}}$.
Alternatively, if one uses
$$B^a\equiv\frac1{V.V}A^a-2\frac{V.A}{(V.V)^2}V^a$$
instead of the na\"{\i}ve acceleration, then (\ref{circles}) becomes
\begin{equation}\label{Pauls_circles}
\textstyle\partial B^a=V.BB^a-\frac12B.BV^a+\Rho_b{}^aV^b\end{equation}
(which is equation (2) from~\cite{T}).
Conformal geodesics were introduced in~\cite{Y} and one can verify, as is done
explicitly in~\cite{BE}, that (\ref{circles})
(equivalently~(\ref{Pauls_circles}) as is done in~\cite{T}), is conformally
invariant. Unlike ordinary geodesics, however, the parameter $\tau$ is
determined only up to projective freedom $\tau\mapsto(a\tau+b)/(c\tau+d)$.

Although dependent on a choice of metric in the conformal class, one can choose
instead to use ordinary arc-length $t:\gamma\to{\mathbb{R}}$ as a parameter.
Besides, as noted in~\cite{C} via the Cartan connection and directly verified
in~\cite{BE}, any smooth curve on a conformal manifold inherits a preferred
family of parameterizations defined up to projective freedom so the location of
conformal geodesics and how they are parameterized can be regarded as separate
issues. As noted in~\cite{T}, an advantage of using arc-length is that it
avoids the spurious feature that a projective parameterization $\tau$ might
blow up for no good reason. Thus, for a choice of metric in the conformal
class, one can insist on the unit vector field $U^a\equiv V^a/\sqrt{V.V}$ along
$\gamma$ with corresponding parameter $t:\gamma\to{\mathbb{R}}$ such that 
\begin{itemize}
\item $\delta t=1$ where $\delta=\partial/\partial t\equiv U^a\nabla_a$
acting on any tensor field,
\item $C^a=\delta U^a$ is the acceleration field along $\gamma$, 
\end{itemize}
and then it is readily verified that
\begin{equation}\label{Pauls_proper_circles}
\delta C^a=\Rho_b{}^aU^b-(C.C+\Rho_{bc}U^bU^c)U^a.\end{equation}
This is the form of the conformal geodesic equation used in \cite{T} to 
conclude that 
\begin{itemize}
\item conformal geodesics always exist for short time,
\item conformal geodesics are uniquely determined by initial conditions:
\begin{itemize}
\item an arbitrary unit velocity vector $U^a|_{t=0}$,
\item an arbitrary acceleration vector $C^a|_{t=0}$ orthogonal to $U^a|_{t=0}$,
\end{itemize}
\item a conformal geodesic can always be continued as long as $C^a$ is finite.
\end{itemize}
Indeed, as observed in~\cite[Theorem~1.1]{T}, the proof is immediate from 
Picard's Theorem and one can also add that, provided the acceleration remains 
bounded and that solutions exist for $0\leq t\leq 1$, 
\begin{itemize}
\item the end point $\gamma^{-1}(1)$ depends smoothly on the initial 
conditions.
\end{itemize}
These conclusions are also valid for conformal geodesics formulated according
to (\ref{circles}) provided that the acceleration but also the preferred
parameter $\tau$ both remain bounded. (Observe that if our chosen metric
$g_{ab}$ is conformally rescaled $\hat{g}_{ab}=\Omega^2g_{ab}$, then our
various notions of acceleration change:
\begin{itemize}
\item $\hat{A}^a=A^a-V.V\Upsilon^b+2(V.\Upsilon)V^b$,
\item $\hat{B}_a=B_a-\Upsilon_a$,
\item $\hat{C}_a=C_a-\Upsilon_a+(U.\Upsilon)U_a$,
\end{itemize}
where $\Upsilon_a=\nabla_a\log\Omega$ but whether any of these
notions remain finite is clearly independent of choice of metric in the
conformal class.)

Notice that if $\gamma\hookrightarrow M$ is a conformal geodesic passing
through $a,b,\in M$ with preferred parameterization
$\tau:\gamma\to{\mathbb{R}}$ so that $\tau(a)=0$ and $\tau(b)=1$, then the same
curve can be also be parameterized by $\hat\tau=1-\tau$. More specifically, the
defining equation (\ref{circles}) holds with $\tau$ replaced by 
$\hat\tau$ but now $\hat\tau(b)=0$ and~$\hat\tau(a)=1$. We shall refer to this 
man{\oe}uvre as {\em reversing the parameterization\/}. 

In Euclidean space, we can solve (\ref{circles}) explicitly. In fact, if one
restricts (\ref{circles}) to a Euclidean plane
${\mathbb{R}}^2\hookrightarrow{\mathbb{R}}^n$, then one may readily verify
that the curve
\begin{equation}\label{probe}\tau\mapsto\big(x(\tau),y(\tau)\big)=
\frac{2}{(2-\alpha\tau)^2+\beta^2\tau^2}
\Big((2-\alpha\tau)\tau,\beta\tau^2\Big)\end{equation}
satisfies 
$$\frac{d}{d\tau} A^a=3\frac{V.A}{V.V}A^a-\frac{3A.A}{2V.V}V^a,$$
where
$$V^a=(dx/d\tau,dy/d\tau)\quad\mbox{and}\quad
A^a=(d^2x/d\tau^2,d^2y/d\tau^2).$$
Hence, this is the unique conformal circle in ${\mathbb{R}}^n$ with initial 
conditions
$$\begin{array}{ll}(x,y,\ldots)|_{\tau=0}=(0,0,\ldots,0)\quad
&(dx/d\tau,dy/d\tau,\ldots)|_{\tau=0}=(1,0,\ldots,0)\\[4pt]
&(d^2x/d\tau^2,d^2y/d\tau^2,\ldots)|_{\tau=0}=(\alpha,\beta,0,\ldots,0)
\end{array}$$
(and one may also verify that for $\beta\not=0$, the trajectory is, indeed, a
round circle in~${\mathbb{R}}^2$ (centred at $(0,1/\beta)$) whilst if
$\beta=0$, the trajectory is the $x$-axis with projective parameterization
$x=2\tau/(2-\alpha\tau)$). Notice that the curve (\ref{probe}) is defined for
all $\tau$ provided $\beta\not=0$, has a limit at
$$\frac{2}{\alpha^2+\beta^2}\Big({-\alpha},\beta\Big)$$
as $\tau\to\pm\infty$, and passes through
$$\frac{2}{(2-\alpha)^2+\beta^2}\Big(2-\alpha,\beta\Big)$$
when $\tau=1$. 

\section{The Uniqueness Theorem}
In this section we prove Theorem~\ref{main_thm}. For convenience, let us denote
by $\Omega$ the open subset $\iota({\mathbb{R}}^n)\subset(M,g)$ that is
supposed conformal to~${\mathbb{R}}^n$. Initially, we know very little about
$\Omega\subset M$:
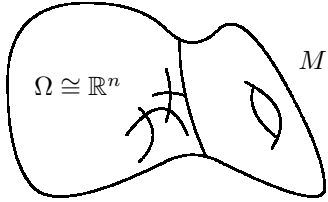
$$\begin{picture}(120,80)\linethickness{.6pt}
\qbezier(0,50) (0,80) (30,80)
\qbezier(30,80) (50,80) (60,70)
\qbezier(60,70) (70,60) (80,70)
\qbezier(80,70) (90,80) (110,40)
\qbezier(110,40) (120,20) (120,10)
\qbezier(120,10) (120,0) (80,20)
\qbezier(80,20) (70,25) (60,20)
\qbezier(60,20) (20,0) (10,10)
\qbezier(10,10) (0,20) (0,50)
\qbezier(90,50) (110,40) (100,25)
\qbezier(93,48) (90,38) (102,29)
\qbezier(65,66) (65,44) (75,22)
\qbezier(67,45) (60,47) (55,45)
\qbezier(60,55) (63,45) (60,35)
\qbezier(45,30) (60,50) (68,30)
\qbezier(50,20) (60,30) (50,40)
\put(10,45){$\Omega\cong{\mathbb{R}}^n$}
\put(110,55){$M$}
\end{picture}$$
In particular, the boundary $\partial\Omega\subset M$ may be very badly
behaved. Nevertheless, by surrounding each point of $\partial\Omega$ by a
co\"ordinate patch $U$ in which one considers Euclidean spheres centred on
points in $U\cap\Omega$ (cf.~\cite[Lemme~14]{Fr}), it follows that in 
$\partial\Omega$ there is a dense set of {\em highly accessible\/} points~$p$,
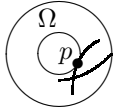
$$\begin{picture}(40,40)\linethickness{.6pt}
\put(20,20){\circle{40}}
\qbezier(40,20) (30,10) (20,10)
\qbezier(25,5) (25,20) (35,25)
\put(20,20){\circle{16}}
\put(24.5,14){$\bullet$}
\put(12,30){$\Omega$}
\put(20,18){$p$}
\end{picture}$$
namely those for which all smooth curves $\gamma(s)$ emanating from $p$ with
velocity on one side of a suitable hyperplane in $T_pM$ remain inside $\Omega$
for $0<s<\epsilon$ for some $\epsilon>0$. In particular, this property is valid
for conformal circles emanating from $p$ with velocity pointing into $\Omega$
in this sense. Let us fix a highly accessible point in $\partial\Omega$ and, at
the risk of putting the cart before the horse, denote it by~$\infty$. We aim to
show that $\partial\Omega=\{\infty\}$. Indeed, if we can show this, then
Theorem~\ref{main_thm} follows easily because then the conformal Weyl
curvature, in case $n\geq 4$, or the Cotton-York tensor, in case $n=3$,
vanishes on $M$ by continuity. This implies that $M$ is conformally flat and
Theorem~\ref{main_thm} follows from Liouville's theorem.

Let us fix a conformal geodesic $\gamma\hookrightarrow M$ with preferred
parameterization $\tau:\gamma\to{\mathbb{R}}$ emanating from $\infty$. By
rescaling the parameterization we may arrange that $\gamma$ is contained inside
$\Omega$ for $0<\tau\leq 1$. Finally, by reversing the parameterization, we
obtain a conformal geodesic starting at a point inside $\Omega$ when $\tau=0$
and ending up at a point $\infty\in\partial\Omega$ when $\tau=1$.

But recall that conformal geodesics and their preferred parameterizations
depend only on the conformal structure and that $\Omega$ is supposed to be
conformal to the standard Euclidean metric on~${\mathbb{R}}^n$. At the end of
the previous section we determined all the conformal circles
in~${\mathbb{R}}^n$ to be straight lines or round circles. A round circle is
bounded and so we are left with having constructed what, from the Euclidean
viewpoint, is simply a straight line. After a Euclidean motion, this straight
line may as well be the $x$-axis with parameterization starting at the origin
and to have the parameter $\tau$ run from $0$ to $1$ at $\infty$ completely
fixes the curve as
\begin{equation}\label{tauoveroneminustau}
\tau\mapsto\Big(\frac{\tau}{1-\tau},0,0,\ldots,0\Big).\end{equation}
In other words, it is the case $(\alpha,\beta)=(2,0)$ in the discussion above. 
Let us compare this curve with the conformal geodesic in 
${\mathbb{R}}^n$ having initial conditions
$$\begin{array}{c}\big(x,y,\ldots\big)|_{\tau=0}=\big(0,0,\ldots,0\big)\qquad
\big(dx/d\tau,dy/d\tau,\ldots\big)|_{\tau=0}=\big(1,0,\ldots,0\big)\\[4pt]
\big(d^2x/d\tau^2,d^2y/d\tau^2,\ldots\big)|_{\tau=0}=
\big(\alpha,\sigma(2-\alpha),0,\ldots,0\big)
\end{array}$$
for some fixed $\sigma$ and for some $0\leq\alpha<2$. According to our previous
discussion, it is the curve (\ref{probe}) with $\beta=\sigma(2-\alpha)$, namely
$$\tau\mapsto
\frac{2}{(2-\alpha\tau)^2+\sigma^2(2-\alpha)^2\tau^2}
\Big((2-\alpha\tau)\tau,\sigma(2-\alpha)\tau^2,0,\ldots,0\Big).$$
When $\tau=1$ this conformal geodesic passes through the point
\begin{equation}\label{endpoint}\frac{2}{2-\alpha}\frac1{1+\sigma^2}
\Big(1,\sigma,0,\ldots,0\Big).
\end{equation} 
Now let $\alpha\uparrow 2$. We see that the endpoint (\ref{endpoint}) moves
monotonely out along the straight line in the direction
$(1,\sigma,0,\ldots,0)$. At least, this is what we see in
$\Omega\cong{\mathbb{R}}^n$. Viewed in $M\supset\Omega$, there is no problem
with the limit curve with initial conditions $(\alpha,\beta)=(2,0)$. By
construction, it is the curve (\ref{tauoveroneminustau}) joining the origin
$0\in{\mathbb{R}}^n\cong\Omega\hookrightarrow M$ to $\infty\in
M\setminus\Omega$. We conclude that the endpoint (\ref{endpoint}) tends to
$\infty\in M$ along the straight line in the direction 
$(1,\sigma,0,\ldots,0)$. At this point, there are several ways to complete the 
proof, one of which is as follows.

Even if we restrict to $|\sigma|\leq 1$ to force uniform convergence
to~$\infty$, we are free to rotate around the $x$-axis to obtain an
entire cone 
$$
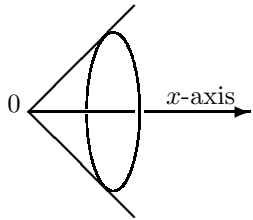
\begin{picture}(100,80)(-50,-40)\linethickness{.6pt}
\qbezier (10., 0.), (10., 12.42), (7.070, 21.21)
\qbezier (7.070, 21.21), (4.14, 30.), (0., 30.)
\qbezier (0., 30.), (-4.14, 30.), (-7.070, 21.21)
\qbezier (-7.070, 21.21), (-10., 12.42), (-10., 1.5)
\qbezier (-10., -1.5), (-10., -12.42), (-7.070, -21.21)
\qbezier (-7.070, -21.21), (-4.14, -30.), (0., -30.)
\qbezier (0., -30.), (4.14, -30.), (7.070, -21.21)
\qbezier (7.070, -21.21), (10., -12.42), (10., 0.)
\thicklines
\put(-32,0){\line(1,0){40}}\put(12,0){\vector(1,0){40}}
\put(-32,0){\line(1,1){40}}
\put(-32,0){\line(1,-1){40}}
\put(20,3){$x$-axis}
\put(-40,0){$0$}
\end{picture}$$
all the rays of which end up at $\infty\in M$. But now we may repeat this
argument starting with any ray from this cone, rather than the $x$-axis. In
this way we obtain a finite collection of cones that cover all rays emanating
from the origin. We conclude that $\infty$ is the only boundary point
of~$\Omega\subset M$.\hfill Q.E.D.

%

\end{document}